\documentclass[12pt]{amsart}
\usepackage{amsmath}%\usepackage{amsrefs}
\usepackage{amsfonts}\usepackage{verbatim}
\usepackage{amssymb}\usepackage{marvosym}
\usepackage{dictsym}\usepackage{wasysym}
\usepackage{amstext}
\usepackage{amsbsy}
\usepackage{amsopn}%\usepackage{MnSymbol}
\usepackage{amsthm}\usepackage{color}

\def\proclaim#1{\vskip0.5em\noindent{\bf #1}\it }
\def\endproclaim{\vskip0.5em\par\noindent\rm}

\def\proclaim#1{\vskip0.5em\noindent{\bf #1}\it}
\def\endproclaim{\vskip0.5em\par\noindent\rm}

\def\undersetbrace#1\to#2{\underbrace{#2}_{#1}}

\def\demo#1{\vskip0.5em\noindent{\bf #1\ }}

\def\text#1{\mbox{#1}}
\def\flushpar{\par\noindent}

\newcommand{\mapright}[1]{%
    \smash{\mathop{%
        \hbox to 1cm{\rightarrowfill}
        }
    \limits^{#1}
    }
}
\newcommand{\mapleft}[1]{%
    \smash{\mathop{%
        \hbox to 1cm{\rightarrowfill}
        }
    \limits_{#1}
    }
}

\def\e{\epsilon}
\def\a{\alpha}
\def\b{\beta}
\def\G{\Gamma}
\def\g{\gamma}
\def\d{\delta}
\def\D{\Delta}
\def\s{\sigma}

\def\th{\theta}
\def\L{\Lambda}
\def\x{\times}

\def \F{\mathcal F}
\def\f{\flushpar}
\def\u{\underline}
\def\v{\varphi}

\def\om{\omega}
\def\Om{\Omega}
\def\B{\mathcal B}
\def\T{\widehat T}
\def\({\biggl(}
\def\){\biggr)}
\def\<{\langle}
\def\>{\rangle}

\def\bul{\smallskip\f$\bullet\ \ \ $}\def\lfl{\lfloor}\def\rfl{\rfloor}

\def\){\biggr)}

\def\<{\bold\langle}
\def\>{\bold\rangle}

\def\bul{\smallskip\f$\bullet\ \ \ $}\def\sms{\smallskip\f}\def\sbul{\f$\bullet\
\ \ $}\def\Lra{\Longrightarrow}
\def\xbm{(X,\B,m)}\def\xbmt{(X,\B,m,T)}\def\xyr{\xrightarrow}\def\isom{\text{\tt Isom}\,}

\begin{document}\title{ On multiple recurrence and other properties of ``nice" infinite measure preserving transformations  }
\author{ Jon. Aaronson $\&$ Hitoshi Nakada}
 \address[Aaronson]{\ \ School of Math. Sciences, Tel Aviv University,
69978 Tel Aviv, Israel.\f\ \ \  {\it Webpage }: {\tt http://www.math.tau.ac.il/$\sim$aaro}}\email{aaro@post.tau.ac.il}
\

\address[Nakada]{ Dept. Math., Keio University,Hiyoshi 3-14-1 Kohoku,
 Yokohama 223, Japan}\email[Nakada]{nakada@math.keio.ac.jp}
\subjclass[2010]{37A40, 53D25, 60F05}
\keywords{Infinite ergodic theory, multiple recurrence,  multiple  rational ergodicity, multiple  rational weak mixing, hyperbolic geodesic flow, Abelian cover.}
\thanks { Aaronson's research was partially supported by ISF grant No. 1599/13. Nakada's research was partially supported by
JSPS grant No. 2430020.\ \ \copyright  2013-5. }

\begin{abstract}
We discuss  multiple versions of rational ergodicity and rational weak mixing for ``nice" transformations, 
including Markov shifts, certain interval maps and hyperbolic geodesic flows.
 These properties entail multiple recurrence.
\end{abstract}
\maketitle\markboth{multiple recurrence and other properties}{Jon. Aaronson $\&$ Hitoshi Nakada}
\section*{\S1. Introduction: multiple properties}
\

The measure preserving transformation ({\tt MPT}) $\xbmt$ is called

 \bul\   \ \ {\tt $d$-recurrent}   as in \cite{A-N}\ \ if $\forall\ A\in\mathcal B,\ m(A)>0\ \exists\ n\ge 1$ so that
$$m(A\cap T^{-n}A\cap T^{-2n}A\cap\dots\cap T^{-dn}A)>0;$$
$1$-recurrence being equivalent to conservativity.

 In \cite{A-N}, the authors considered the multiple recurrence of the
 {\bf Markov shift}  $\xbmt$  of the stochastic matrix $P\colon S\x S\to [0,1]$  with
invariant distribution $\{\mu_s\colon s\in S\}$  where
\begin{align*}X=S^{{\Bbb Z}},\ T&=\text{\tt shift}, \B=\s(\{\text{\tt cylinders}\})\  \& \\ &
m([s_0,\dots,s_n]_k)=\mu_{s_0}p_{s_0,s_1}\dots p_{s_{n-1},s_n}\\ & \text{where the {\it cylinder}}\ [s_0,\dots,s_n]_k:=\{x\in X\colon x_{k+j}=s_j\ \forall 0\le j\le n\},
\end{align*}
showing for $d\in\Bbb N$ that if $\xbmt$ is conservative, ergodic, then
\bul $T$ is
$d$-recurrent $\Leftrightarrow$
$\underset{\text{\tiny $d$-times}}{\underbrace{T\x\dots\x T}}$ is
conservative, ergodic.
\

For Markov shifts, $d$-recurrence is equivalent to {\tt $d$-rational ergodicity}:
\

The conservative, ergodic, measure preserving transformation ({\tt CEMPT}) $\xbmt$ is called {\it  $d$-rationally ergodic   along $\frak K\subset\Bbb N$} if there exist:
 \sms (i) a sequence of constants $a_d(n)\uparrow\ \infty$ and
 \sms (ii) a dense, $T$-invariant, hereditary ring (aka ideal)  $R_d(T)\subset\mathcal F:=\{F\in\B:\
 m(F)<\infty\}$  s.t.
 \begin{align*}\frac1{a_d(n)}\sum_{k=0}^{n-1}m(\bigcap_{j=0}^dT^{-(jk+r_j)}B_j)&\xyr[n\to\infty,\ n\in\frak K]{}\prod_{j=0}^dm(B_j)\\ & \forall\ B_0,B_1,\dots,B_d\in R_d(T)\ \&\ r_0,\dots, r_d\in\Bbb Z;
 \end{align*}
 
For Markov shifts,  
$$u_d(n)=u_n^d,\ \ a_{d}(n) = \sum_{k=1}^{n} u_{n}^{d}$$ where $u_n:=\frac{p_{s,s}^{(n)}}{\mu_s}$ (any $s\in S)$ and
$$R_d(T)\supset\{B\in\mathcal F:\ B\subset C \ \text{\tt a cylinder}\}.$$
In this paper we extend this to further classes of  transformations which we call  ``{\tt nice}''.
\

 Some notes on terminology:
 \bul $d$-{\it rational ergodicity} is  $d$-rational ergodicity along $\Bbb N$;
 \bul {\it subsequence $d$-rational ergodicity} is  $d$-rational ergodicity along some $\frak K\subset\Bbb N$;
 \bul $1$-rational ergodicity  is called {\it weak rational ergodicity} in \cite{RatErg}.
\

Evidently, subsequence $d$-rational ergodicity implies $d$-recurrence

\

\subsection*{Description of results}
\

In \S2, we define ``$d$-nice transformations'' and show that these are subsequence $d$-rational ergodic. This is applied in 
\S5 where we establish $2$-recurrence [{\tt $2$-dissipation}] of
the geodesic flow of a $\Bbb Z$-cover [{\tt $\Bbb Z^2$-cover}] 
of a compact, hyperbolic surface, as advertised in \cite{PPS}.
In \S3, we give sufficient conditions for the stronger  multiple rational weak mixing properties which are applied to show
 $1$-rational weak mixing 
 \sms - of certain special semiflows  in 
 \S4; 
 \sms - of the geodesic flow of  a $\Bbb Z^\kappa$-cover ($\kappa=1,2$)
of a compact, hyperbolic surface.

\

\subsection*{Asymptotic equivalence of sequences}\ \ Throughout the paper we use the following notations for sequences 
$(a_1,a_2,\dots),\ (b_1,b_2,\dots)\in\Bbb R^\Bbb N$:
\bul $a_n\ \approx\ b_n$ if $a_n-b_n\xrightarrow[n\to\infty]{}\ 0$,
\bul for $M>0,\ a_n=b_n\pm M$ if $|a_n-b_n|\le M\ \forall\ n\ge 1$;
\

and for $(a_1,a_2,\dots),\ (b_1,b_2,\dots)\in\Bbb R_+^\Bbb N$:
\bul $a_n\ \sim\  b_n$ if $\frac{a_n}{b_n}\xrightarrow[n\to\infty]{}\ 1$,
\bul $a_n\ \propto\  b_n$ if $\frac{a_n}{b_n}\xrightarrow[n\to\infty]{}\ c\in\Bbb R_+$,
\bul $a_n\ll b_n$ if $\exists\ M>0$ so that $a_n\le Mb_n\ \forall\ n\ge 1$,
\bul $a_n\asymp b_n$ if $a_n\ll b_n$ and $b_n\ll a_n$.
\

For $M>1$ we also write
\bul $a_n=M^{\pm 1}b_n$ or $a_n\overset{M}\asymp b_n$ if $\frac{b_n}M\le a_n\le Mb_n\ \forall\ n\ge 1$.
\section*{\S2 Nice transformations}
\

Let $\xbmt$  be a {\tt CEMPT}.  %For $d\ge 1$, w
\

We'll call  a set $\Om\in\mathcal F_+:=\{A\in\B:\ 0<m(A)<\infty\}$ {\it admissible} for $T$ if
\sms   $m(\bigcap_{k=0}^dT^{-kn}\Om)\asymp u(n)^d\ \ \forall\ d\ge 1$ where
$u(n)=u(\Om,n):=\tfrac{m(\Om\cap T^{-n}\Om)}{m(\Om)}$.

\

Let  $d\in\Bbb N$. We'll call the {\tt CEMPT} $\xbmt$ {\it  $d$-nice} if 
\sms (i) there is an admissible set $\Om\in\mathcal F_+$  for $T$;
$$\text{(ii)}\ \ (X_d,\B_d,m_d,T_d):=(X^d,\B(X^d),
\underset{\text{\tiny $d$-times}}{\underbrace{m\x\dots\x m}},\underset{\text{\tiny $d$-times}}{\underbrace{T\x\dots\x T}})\ \ \text{is a {\tt CEMPT};}$$
\sms (iii) $\exists\ M>1$ and a countable, dense collection $\mathcal A\subset\B\cap \Om$ with $\Om\in\mathcal A$ s.t. $\forall\ r_0,\dots,r_d\in\Bbb Z$,
$$\sum_{k=0}^{n-1}m(\bigcap_{k=0}^dT^{-kn+r_k}B_k)\overset{M}{\asymp}\prod_{k=0}^dm(B_k)a_d(n)\ \forall\ B_0,B_1,\dots,B_d\in\mathcal A$$
where $a_d(n):=\sum_{k=1}^n u(\Om,k)^d\to\infty$ by  condition  (ii), 
\

\sms (iv)  $a_d(2n)\ll a_d(n)$.
\

\proclaim{Proposition 2.1}\ \  If is $T$ $d$-nice, then $a_d(n)\to\infty$ and
$T_d$ is $1$-rationally ergodic with return sequence $a_n(T_d)\asymp a_d(n)$.\endproclaim
\demo{Proof}\ \ For $\Om\in\F_+$ admissible,
$$\int_{\Om^d}S_n^{(T_d)}(1_{\Om^d})^2dm_d\ll \(\sum_{k=1}^nu(\Om,k)^d\)^2\asymp \(\int_{\Om^d}S_n^{(T_d)}(1_{\Om^d})dm_d\)^2.\ \ \ \CheckedBox$$
\proclaim{Theorem 2.2}
\

 \ If $\xbmt$ is $d$-nice, then it is subsequence $d$-rationally ergodic.
and the hereditary ring satisfies
$$R_d(T)\supset\{B\in\mathcal F:\ \exists\ n\in\Bbb Z,\ B\subset T^{n}\Om\}.$$\endproclaim\
\

The proof of theorem 2.2 uses:

\proclaim{ Lemma 2.3}\ \ Let $\xbmt$ be $d$-nice and let $\Om\in\F_+$ be admissible. For any $0\le \nu\le d$, define
$$\psi^{(\nu)}_n:=\sum_{k=1}^n\prod_{i=1}^\nu 1_\Om\circ T^{-ik}\cdot
\prod_{j=1}^{d-\nu} 1_\Om\circ T^{jk},$$
then
$$\text{\rm (i)}\ \ \int_\Om\psi^{(\nu)}_ndm\asymp\ a_d(n)\ \ \&\ \ \text{\rm (ii)}\ \ \int_\Om(\psi^{(\nu)})^2dm =O\(a_d(n)^2\).$$\endproclaim

 \subsection*{ Proof of lemma 2.3} \ This lemma is a generalization of lemma 1.5 in \cite{A-N}.

Throughout, we use
admissibility of $\Om$:
\f if  $b(1),\dots,b(\kappa)\in\Bbb Z$ and $b(1)\le b(2)\le\dots\le b(\kappa)$
then
\begin{align*}
 \tag0 m\(\bigcap_{r=1}^{\kappa}T^{-b(r)}\Om\)=m\(\bigcap_{r=1}^{\kappa}T^{b(r)}\Om\)\asymp
\prod_{r=2}^{\kappa }m(\Om\cap T^{-(b(r)-b(r-1))}\Om).
\end{align*}

Set
$$\e_k(\nu):=\prod_{-\nu\le j\le d-\nu,\ j\ne 0}
1_\Om\circ T^{jk},$$
then
$$\psi^{(\nu)}_n
=\sum_{k=1}^n\e_k(\nu)\ \text {and }\ 
\int_\Om(\psi^{(\nu)}_n)^2dm\le 2\sum_{k=1}^n\sum_{\ell=k}^n
\int_\Om\e_k(\nu)\e_\ell(\nu)dm.$$

In view of (0), the form of $\int_\Om\e_k(\nu)\e_\ell(\nu)dm$ depends on the orders of
the sets $\{ik,j\ell\colon \ 1\le i,j\le\nu\}$ and
$\{ik,j\ell\colon \ 1\le i,j\le d-\nu\}$.

The ordering of 
$$\Om_d(k,\ell):=\{ik,\ j\ell\colon \ 1\le i,j\le d\}\ \ \ (d,k,\ell\in\Bbb N)$$ is treated in lemma 1.5 of \cite{A-N} and we'll use
results established there.
\

Define  $N_{(k,\ell)}\colon \Bbb N\x\{0,1\}\to\Bbb N$ by $N_{k,\ell}(j,\e)=
(1-\e)jk+\e j\ell$, then for each $d\ge 1,\  N_{(k,\ell)}\colon \Bbb N_d\x\{0,1\}\to\Om_d(k,\ell)$ is a surjection. Here
$\Bbb N_d:=\{1,2,\dots,d\}$.
\

An {\it ordering} of $\Om_d(k,\ell)$ 
is a bijection $\pi:\Bbb N_{2d}\to \Bbb N_d\x\{0,1\}$ so that $N_{(k,\ell)}\circ\pi :\Bbb N_{2d}\to \Om_d(k,\ell)$ is non-decreasing.
\

Given a bijection $\pi:\Bbb N_{2d}\to \Bbb N_d\x\{0,1\}$, let 
$$D(\pi):=\{(k,\ell)\in\Bbb N\x\Bbb N\  \pi \ \text{\rm orders}\ \Om_d(k,\ell)\}.$$
Possibly $D(\pi)=\emptyset$.  However, as shown in lemma 1.5 in \cite{A-N}, if
$$F_d:=\{{p\over q}\colon \ 0\le p\le q\le d,\ (p,q)=1\}=\{0:=r_0^{(d)}<r_1^{(d)}<\dots<r_{N_d}^{(d)}=1\}$$ is the
{\bf Farey sequence} of order $d$, then for each $0\le j<N_d$, there is a bijection $\pi_j:\Bbb N_{2d}\to \Bbb N_d\x\{0,1\}$ so that
\begin{align*}
\tag1 D(\pi_j)=\{(k,\ell)\in\Bbb N^2\colon \
{k\over\ell}\in (r_j,r_{j+1}]\}.
\end{align*}

\bigskip
\par Let $\frak o_d:=\{\pi_j:\ 0\le j<N_d\}$. It follows from (1) that if $1\le d'<d$ and
$\pi\in\frak o_d,$ then $\exists\ \pi'\in\frak o_{d'}$ such that
$D(\pi)\subset D(\pi')$.
\par Fix $\pi\in\frak o_d,\ (k,\ell)\in D(\pi)$, then  writing $\pi(j)=(\kappa_j,\e_j)$ for $1\le j\le 2d$, we have (see \cite{A-N})
\begin{align*} 
N_{(k,\ell)}\circ\pi(j)-N_{(k,\ell)}\circ\pi(j-1) &=
\kappa_j\big[(1-\e_j)k+\e_j\ell\big]-
\kappa_{j-1}\big[(1-\e_{j-1})k+\e_{j-1}\ell\big]\\ &=\<a_j,(k,\ell)\>\end{align*}
where $N_{(k,\ell)}\circ\pi(0):=0,\ a_1=(\kappa_1(1-\e_1),\kappa_1\e_1)$ and
$$a_j=a_j(\pi):=(\kappa_j(1-\e_j)-\kappa_{j-1}(1-\e_{j-1}),
\kappa_j\e_j-\kappa_{j-1}\e_{j-1})\ \ \ \ (j\ge 2).$$

%\newpage
\par The vectors $\{a_j(\pi)\}_{j=1}^{2d}$ are non-zero and 
$$\{a_j(\pi)\colon \ 1\le j\le 2d\}=\{a_j^{(1)}(\pi),a_j^{(2)}(\pi)\colon \ 1\le j\le d\}$$
where $a_j^{(1)}(\pi)$ and $a_j^{(2)}(\pi)$ are linearly independent
$\forall\ 1\le j\le d\}$.

\

Now, set
$$\e_k^{\pm}(\nu)=\prod_{j=1}^\nu 1_\Om\circ T^{\pm jk},$$
then $\e_k(\nu)=\e_k^-(\nu)\e_k^+(d-\nu)$, 
$$\int_\Om\e_k(\nu)\e_\ell(\nu)dm =
\int_\Om (e_k^-(\nu)\e_\ell^-(\nu))(\e_k^+(d-\nu)\e_\ell^+(d-\nu))dm$$
and it follows from admissibility that
\begin{align*} \int_\Om (e_k^-(\nu)\e_\ell^-(\nu)) &(\e_k^+(d-\nu)\e_\ell^+(d-\nu))dm\ \asymp
\\ &\int_\Om\e_k^+(\nu)\e_\ell^+(\nu)dm\int_\Om\e_k^+(d-\nu)\e_\ell^+(d-\nu)dm.
\end{align*}

\par Using the discussion above (and admissibility of $\Om$), 
\begin{align*}
\tag2 \int_\Om\e_k^+(d)\e_\ell^+(d)dm\ \asymp\ 
\prod_{j=1}^{2d}u(\<a_j,(k,\ell)\>)\ \forall\ (k,\ell)\in D(\pi),\
\pi\in\frak o_d. 
\end{align*}

\

\par We now complete the
proof of lemma 2.3 in case 
 $\nu\ge d-\nu$ (the other case being similar). 
\

For each $\pi\in\frak o_\nu$, let
$\pi'\in\frak o_{d-\nu}$ be such that $D(\pi)\subset D(\pi')$.
\par Since $\{(k,\ell)\in\Bbb N^2\colon \ k\le\ell\}=
\biguplus_{\pi\in\frak o_\nu}D(\pi)$, we have:
\begin{align*} \int_\Om(\psi^{(\nu)}_n)^2dm &\le 2\sum_{k=1}^n\sum_{\ell=k}^n
\int_\Om\e_k(\nu)\e_\ell(\nu)dm\\ &=
2\sum_{\pi\in\frak o_\nu}\sum_{(k,\ell)\in D(\pi),\ k,\ell\le n}
\int_\Om\e_k(\nu)\e_\ell(\nu)dm.\end{align*}

\par For each $\pi\in\frak o_\nu$,
\begin{align*} &\sum_{(k,\ell)\in D(\pi)\cap\Bbb N_n^2}\int_\Om\e_k(\nu)\e_\ell(\nu)dm \\ &=
\sum_{(k,\ell)\in D(\pi)\cap\Bbb N_n^2}
\int_\Om\e_k^+(\nu)\e_\ell^+(\nu)dm
\int_\Om\e_k^+(d-\nu)\e_\ell^+(d-\nu)dm\\ &\asymp
\sum_{(k,\ell)\in D(\pi)\cap\Bbb N_n^2}
\prod_{j=1}^{2\nu}u(\<a_j(\pi),(k,\ell)\>)\prod_{j=1}^{2(d-\nu)}
u(\<a_j(\pi'),(k,\ell)\>)
\\ &=
\sum_{(k,\ell)\in D(\pi)\cap\Bbb N_n^2}
\prod_{j=1}^{d}u(\<a_j^{(1)},(k,\ell)\>)u(\<a_j^{(2)},(k,\ell)\>)\end{align*}
 where
$$\{a_j^{(1)},a_j^{(2)}\}_{j=1}^{2d}=
\{a_j^{(1)}(\pi),a_j^{(2)}(\pi)\}_{j=1}^{2\nu}\cup
\{a_j^{(1)}(\pi'),a_j^{(2)}(\pi')\}_{j=1}^{2(d-\nu)}.$$
\par Consider $B_j\colon \Bbb R^2\to\Bbb R^2$ defined by
$(B_jx)_i:=\<x,a_j^{(i)}\>\ \ \ (i=1,2)$ which is injective.
Let $K>0$ be such that $\|B_jx\|_\infty\le K\|x\|_\infty\ \forall\ x,j$.
\par  By H\"older's inequality,
\begin{align*} & \sum_{(k,\ell)\in D(\pi)\cap\Bbb N_n^2}\prod_{j=1}^{d}u(\<a_j^{(1)},(k,\ell)\>)u(\<a_j^{(2)},(k,\ell)\>)
\\ &\le
\prod_{j=1}^{d}
\(\sum_{(k,\ell)\in D(\pi)\cap\Bbb N_n^2}
u(\<a_j^{(1)},(k,\ell)\>)^du(\<a_j^{(2)},(k,\ell)\>)^d\)^{{1\over d}}
\\ &=
\prod_{j=1}^{d}
\(\sum_{(k,\ell)\in B_j(D(\pi)\cap\Bbb N_n^2)}
u(k)^du(\ell)^d\)^{{1\over d}}\\ &\le
\sum_{(k,\ell)\in \Bbb N_{Kn}^2}u(k)^du(\ell)^d\\ &=a_d(Kn)^2\\ &\ll a_d(n)^2.\ \ \CheckedBox\end{align*}
\demo{Proof of theorem 2.2}
For any $0\le \nu\le d$, fix $A_0,\dots,A_d,B_0,\dots,B_d\in\mathcal A$ so that
$A_j=B_j\ \forall\ 0\le j\le d,\ j\ne\nu$, then
\begin{align*}
 \sum_{k=0}^n|m(\bigcap_{j=1}^dT^{-jk}A_j)-m(\bigcap_{j=1}^d&T^{-jk}B_j)| \\ &=
 \sum_{k=0}^nm(\bigcap_{j=1,\ j\ne\nu}^dT^{-jk}A_j\cap T^{-\nu k}(A_\nu\D B_\nu))\\ &
 \le\int_{A_\nu\D B_\nu}\psi^{(\nu)}_ndm\\ &\le\sqrt{m(A_\nu\D B_\nu)}\sqrt{\int_\Om(\psi^{(\nu)}_n)^2dm}\\ &\le
 M\sqrt{m(A_\nu\D B_\nu)}a_d(n).\ \ \ \CheckedBox
\end{align*}

\section*{\S3 Multiple rational weak mixing.}
\

Let $d\in\Bbb N$. We'll call the {\tt CEMPT} $\xbmt$  {\it  $d$-rationally weakly mixing along $\frak K\subset\Bbb N$}  
if it is $d$-rationally ergodic and
\bul $\exists$ $u_d(n)>0$ so that the normalizing constants are given by $a_d(n):=\sum_{k=0}^{n-1}u_d(n)\uparrow\ \infty$ where
  \begin{align*}\frac1{a_d(n)}\sum_{k=0}^{n-1}&\left|m(\bigcap_{j=0}^dT^{-(jk+r_j)}B_j)-\prod_{j=0}^dm(B_j)u_d(k)\right| 
  \xyr[n\to\infty,\ n\in\frak K]{}\ 0\end{align*}

  \
  
  Rational weak mixing (i.e. $1$-rational weak mixing along $\Bbb N$) was introduced  in \cite{RatWM}.
  \
  
  Recall that the {\tt CEMPT} $\xbmt$ is {\it pointwise dual ergodic} if there are constants $a_n(T)>0$ so that
 \begin{align*}
 \tag{ DK}\ \ \frac1{a_n(T)}\sum_{k=1}^n\widehat{T}^kf\xyr[n\to\infty]{}\int_Xfdm\ \ \text{a.e.}\ \ \forall\ f\in L^1(m).
\end{align*} 
Here, $\widehat T:L^1(m)\to L^1(m)$ is
 the {\it transfer operator}  defined by  $\widehat Tf:=\frac{d\nu_f\circ T^{-1}}{d\,m}$ where $\nu_f(A):=\int_Xfdm$. It satisfies
$$\int_X \widehat Tf\cdot gdm=\int_X f\cdot  Tgdm\ \forall\ f\in L^1(m),\ g\in L^\infty(m).$$
\subsection*{First return time and induced transformation}
\

Suppose $(X,\B,m,T)$ is a {\tt CEMPT} and let
$\Om\in \B,\ m(\Om)>0$, then  $m$--a.e. point of $\Om$ returns
 to $\Om$ under iterations of $T$.
 The {\it return time} function to $\Om$, defined for
$x\in \Om$ by $\v_\Om(x):= \min\{n\ge 1: T^nx\in \Om\}$ is finite
$m$--a.e. on $\Om$.
\par The {\it induced transformation} on $\Om$ is defined by
$T_\Om x= T^{\v_\Om(x)}x$.
\par In case $m(\Om)<\infty$, then  $(\Om,\B\cap \Om,T_\Om,m_\Om)$ is an {\tt EPPT} (ergodic, probability preserving transformation).
 
\proclaim{Proposition 3.1}\ \  

Let  $\xbmt$ be an exact, pointwise dual ergodic {\tt CEMPT} and suppose that $\exists\ \Om\in\mathcal F_+\ \&\ \a\subset\B\cap\Om$ a one-sided $T_\Om$-generator for $\B\cap\Om$ so that
\sms (i) the  first return time $\v_\Om:\Om\to\Bbb N$ is $\a$-measurable $\&$
\sms (ii)  there exist\ $d\in\Bbb N,\  \ 0<\g <\tfrac1d$ and $b:\Bbb R_+\to\Bbb R_+$, a $\tfrac1\g $-regularly varying function, so that 
\begin{align*}\text{\tt(LLT)}\ \ 
 b(n)\widehat{T}_\Om^n(1_{A\cap[\v_n=k_n]})&\xyr[n\to\infty,\ \frac{k_n}{b(n)}\to x]{L^\infty(m_\Om)}m_\Om(A)f_\g (x)\\ & \forall\ \text{\it cylinders}\  A\ \&\  x\in\Bbb R_+,
\end{align*}
where $f_\g $ is the probability density function of the normalized, positive $\g $-stable random variable, then
\begin{align*}
\tag{\dsagricultural} \frac1{a_d(n)}\sum_{k=1}^n|\widehat{T}^{k+r_1}(1_{A_1}&\widehat{T}^{k+r_2}(1_{A_2}\dots\widehat{T}^{k+r_d}(1_{A_d})\dots)
-\prod_{i=1}^dm(A_i)u_k^d|\\ &\xyr[n\to\infty]{}0\ \ \text{a.e.} \ \ \forall\ A_1,\dots,A_d\in\B.
\end{align*}  and, in particular, $\xbmt$ is $d$-rationally weakly mixing.\endproclaim
\subsection*{Remark}\ \ \ Suitable:
\sbul {\tt AFN} maps, and towers over {\tt AFU maps} as in \cite{ADSZ},\sbul  towers over {\tt Gibbs-Markov maps} as in \cite{AD} 
\

satisfy the assumptions of proposition 3.1.
\

\demo{Proof}
\

Let $a(n):=b^{-1}(n),\ u_n:=\frac{\g  a(n)}n$, then
$$a_n(T)\sim a(n)\sim\sum_{k=1}^n u_k.$$

We'll use pointwise dual ergodicity  to establish rational weak mixing.

To this end, we show first (as in \cite{G-L}) that ({\tt LLT}) implies that :
 $$\text{(\Football)}\ \ \varliminf_{n\to\infty}\frac1{u_n}\widehat{T}^n1_A\ \ge m(A)\ \ \text{a.e. for $A\subset\Om$ a union of cylinders}.$$
\demo{Proof of (\Football)}
\

As in \cite{G-L} (see also \cite{RatWM}):
\begin{align*} \ \T^{n}1_A&=\sum_{k=1}^n\T_\Om^{k}1_{A\cap[\v_k=n]}\\ &\ge
\sum_{1\le k\le n,\ x_{k,n}\in [c,d]}\T_\Om^{k}1_{A\cap[\v_k=x_{k,n}{b(k)}]}\ \ \text{where}\ x_{k,n}:=\tfrac{n}{b(k)}
\\ &\overset{\text{\small(LLT)}}{\text{\Large $\sim$}}
\sum_{1\le k\le n,\ x_{k,n}\in [c,d]}\tfrac{f(x_{k,n})}{b(k)}m(A).
\end{align*}
By the $\g$-regular variation of $a$,
$$\frac1{b(k)}\ \text{\Large $\sim$}\ \ \frac{\g a(n)}n\cdot\frac{x_{k,n}-x_{k+1,n}}{x_{k,n}^\g}.$$
Therefore\ \ 
\begin{align*}\frac1{u_n}\sum_{1\le k\le n,\ x_{k,n}\in [c,d]}\frac{f(x_{k,n})}{b(k)}\ \ &\text{\Large$\approx$}
\sum_{1\le k\le n,\ x_{k,n}\in (c,d)}\frac{(x_{k,n}-x_{k+1,n})}{x_{k,n}^\g}f(x_{k,n})\\ &\xrightarrow[n\to\infty]{}\ \ 
\int_{[c,d]}\frac{f(x)dx}{x^\g}\\ &=\ \ \mathbb E(1_{[c,d]}(Z_\g)Z_\g^{-\g})\\ &\xrightarrow[c\to 0+,\ d\to\infty]{}\ \ 1.\ \ \CheckedBox\text{\Football}\end{align*}
\demo{Proof of\ \ (\dsagricultural)}\ \ \

Let  $A\subset\Om$ be a finite unions of cylinders. By (DK), (\Football) and proposition 3.3 in \cite{RatWM},
$$\frac1{a(n)}\sum_{k=1}^n|\widehat{T}^{n+r}1_{A}-m(A)u_n|
\xyr[n\to\infty]{}0\ \ \text{a.e.}$$
whence by proposition 3.1 in \cite{RatWM}, for a.e. $x\in\Om$, there is a subset $K=K_x\subset\Bbb N$ of full density so that
$$\frac1{u_n}\widehat{T}^{n+r}1_{A}(x)
\xyr[n\to\infty,\ n\in K_x]{}m(A).$$
 Let $C\subset \Om$ be compact.
As in \cite{RatWM}, $\exists$ finite unions of cylinders $A_n\supset A_{n+1}\supset C$ so that $m(A_n)\downarrow m(C)$ as $n\to\infty$. Thus,
for a.e. $x\in\Om$, there is a subset $K=K_x\subset\Bbb N$ of full density so that
$$\varlimsup_{n\to\infty,\ n\in K_x}\frac1{u_n}\widehat{T}^{n+r}1_{C}(x)\le m(C),$$
whence, again by (DK)  and propositions 3.1 $\&$ 3.3 in \cite{RatWM},
$$\frac1{u_n}\widehat{T}^{n+r}1_{C}(x)
\xyr[n\to\infty,\ n\in K_x]{}m(C).$$
 Let\ $B\in\B(\Om)$. Again as in \cite{RatWM}, there are  compact sets $C_n\subset C_{n+1}\subset B$ so that $m(C_n)\uparrow m(B)$ as $n\to\infty$. 
Thus, for a.e. $x\in\Om$, 
there is a subset $K=K_x\subset\Bbb N$ of full density so that
$$\varliminf_{n\to\infty,\ n\in K_x}\frac1{u_n}\widehat{T}^{n+r}1_{B}(x)\ge m(B),$$
whence, (as above) for a possibly smaller $K=K_x\subset\Bbb N$ of full density,
$$\frac1{u_n}\widehat{T}^{n+r}1_{B}(x)
\xyr[n\to\infty,\ n\in K_x]{}m(B).$$

\

It follows that for   $A_1,\dots,A_d\in\B(\Om)$, for a.e. $x\in\Om$, there is a subset $K=K_x\subset\Bbb N$ of full density so that
$$\frac1{u_n^d}\widehat{T}^{k+r_1}(1_{A_1}\widehat{T}^{k+r_2}(1_{A_2}\dots\widehat{T}^{k+r_d}(1_{A_d})\dots)
\xyr[n\to\infty,\ n\in K_x]{}\prod_{i=1}^dm(A_i).$$
The index of regular variation of $u_n^{-d}$ is $d\g\in (0,1)$ so, again by  propositions 3.1 and 3.3 in \cite{RatWM},
\begin{align*}
\frac1{a_d(n)}\sum_{k=1}^n|\widehat{T}^{k+r_1}(1_{A_1}&\widehat{T}^{k+r_2}(1_{A_2}\dots\widehat{T}^{k+r_d}(1_{A_d})\dots)
-\prod_{i=1}^dm(A_i)u_k^d|\\ &\xyr[n\to\infty]{}0\ \ \text{a.e.}. \ \ \ \CheckedBox
\end{align*}

\section*{\S4  Special semiflows}
\

Given a $\s$-finite measure space $\xbm$ we denote
$$\text{\tt MPT}\xbm:=\{\text{measure preserving transformations of}\ \xbm\}.$$
Note that these transformations are not necessarily invertible.
In this section, we consider, for $S$ a finite set and $\kappa\in\Bbb N$, measure preserving semiflows  $\Psi:\Bbb R_+\to\text{\tt MPT}\xbm$  where
\begin{align*}
 &X=\{(x,n,t)\in\Om\x\Bbb Z^\kappa\x\Bbb R_+:\ 0\le x<h(x)\}\ \ \ \text{where}\\ &
 \Om\subset S^\Bbb N\ \ \text{a transitive {\tt SFT},\ $h:\Om\to\Bbb R_+$  H\"older};\\ & \B=\B(X),\ \ \ m(A\x B\x C)=\mu(A)\#(B)\text{\tt Leb}\,(C)
 \\ &\text{where $\#$ is counting measure,  $\mu\in\mathcal P(\Om)$ is Gibbs as in \cite{Rees}, \cite{BowenBook};}\\ &
 \Psi_t(x,n,y)=(T^nx,n+\phi_n(x),y+t-h_n(x))\\ &\  \text{where $\phi:\Om\to\Bbb Z^\kappa$ is continuous, $T=\text{\tt Shift}\  \& \ n=n_t(x,y)$ is s.t.}\\ & h_n(x):=\sum_{k=0}^{n-1}h(T^kx)\le y+t<h_{n+1}(x).
\end{align*}

\subsection*{Pointwise dual ergodicity}

Suppose that, $\Psi$ is ergodic, equivalently $T_\phi$ is ergodic, or $\phi$ is {\it non-arithmetic} in the sense that
$\nexists$ solution to
\begin{align*}&\phi=k+g-g\circ T,\ \ g:\Om\to\Bbb Z^\kappa\ \ \text{measurable},\\ & k:\Om\to \Bbb K\ \ 
\text{measurable, where $\Bbb K$ is a proper subgroup of $\Bbb Z^\kappa$}. 
\end{align*}

By the central limit theorem  in \cite{G-H}, 
\begin{align*}\tag{{\tt CLT}} \frac{\phi_n}{\sqrt n}\xrightarrow[n\to\infty]{\text{\tiny distribution}}\ X
\end{align*}
where $X$ is a globally supported, centered Gaussian random variable on $\Bbb R^\kappa$.
\

As in \cite{mode},  $T_\phi$ is dissipative for $\kappa\ge 3$ and
pointwise dual  ergodic for $\kappa=1,2$ with

\begin{align*}
a_n(T_\phi)\sim\sum_{k=0}^nu_k(T_\phi)\ \sim \begin{cases} &2\sqrt nf_X(0)\ \ \ \ \ \ \ \kappa=1;\\ &\log nf_X(0)\ \ \ \ \ \ \ \kappa=2.\end{cases}
\end{align*}
Analogously to
Proposition 2.2 in \cite{mode}, 
$(X,\B,m,\Psi)$ is dissipative for $\kappa\ge 3$ and
pointwise dual  ergodic (as a flow) for $\kappa=1,2$ with
$$a_n(\Psi)\sim \varkappa^{\frac{\kappa}2-1}a_{n}(T_\phi)$$
where $\varkappa=\int_\Om hdP\in\Bbb R_+$: namely
\begin{align*}
    \frac1{a_n(\Psi)}\int_0^n\widehat{\Psi}_t(F)dt\xrightarrow[n\to\infty]{\text{a.e.}}\ \int_X Fdm \ \forall\ F\in L^1(m).
\end{align*}
\proclaim{Proposition 4.1 \ \ (Exactness)}\ \ Suppose that the function $(h,\phi):\Om\to\Bbb G:=\Bbb R\x\Bbb Z^\kappa$ is  {\tt  non-arithmetic} in the sense that
$\nexists$ solution to
\begin{align*}&(h,\phi)=k+g-g\circ T,\ \ g:\Om\to\Bbb G\ \ \text{measurable},\\ & k:\Om\to \Bbb K\ \ 
\text{measurable, where $\Bbb K$ is a proper subgroup of $\Bbb G$}, 
\end{align*}
then $(X,\B,m,\Psi_t)$ is an exact endomorphism $\forall\ t>0$.\endproclaim\demo{Proof} \ \ The assumption is equivalent to the ergodicity of 
$$(\Om\x\Bbb G,\mu\x\text{\tt Leb}\x\#,T_{(h,\phi)})$$ which entails (characterizes) the exactness of $\Psi$.\ \ \Checkedbox
\

In this case, 
 for each $t>0$,
$(X,\B,m,\Psi_t)$ is 
pointwise dual  ergodic (as a transformation) for $\kappa=1,2$ with
$$a_n(\Psi_t)\sim \varkappa^{\frac{\kappa}2-1}a_{tn}(T_\phi),$$
 namely
\begin{align*}\tag{\tt PDE}
    \frac1{a_n(\Psi_t)}\sum_{k=0}^n\widehat{\Psi}_{tk}(F)\xrightarrow[n\to\infty]{\text{a.e.}}\ \int_X Fdm \ \forall\ F\in L^1(m).
\end{align*}
\subsection*{Rational weak mixing}
\

If   $\phi$ is {\it aperiodic} in the sense that
$\nexists$ solution to
\begin{align*}&\phi=k+g-g\circ T,\ \ g:\Om\to\Bbb Z^\kappa\ \ \text{measurable},\\ & k:\Om\to \Bbb K\ \ 
\text{measurable, where $\Bbb K$ is a proper coset of $\Bbb Z^\kappa$}, 
\end{align*}
then by theorem C* in \cite{HH} (which in this case follows from the local limit theorem in \cite{G-H}), for $A\subset\Om$ a cylinder set:
\begin{align*}\tag{{\tt LLT}} \T^n(1_{A\cap [\phi_n=\lfl t_n\rfl]})\sim\ \frac{f_X(\frac{t_n}{\sqrt n})\mu(A)}{n^{\frac{\kappa}2}} \ \ 
\text{as}\ n\to\infty\ \&\ t_n=O(\sqrt n).
\end{align*}
In particular,
$$\T_\phi(1_{A\x\{0\}})=\T^n(1_{A\cap [\phi_n=0]})\sim\ \frac{\mu(A)f_X(0)}{n^{\frac{\kappa}2}}=:u_n(T_\phi)\mu(A) \ \ \text{as}\ n\to\infty.$$
It follows from proposition 3.1 that $(\Om\x\Bbb Z^\kappa,\B(\Om\x\Bbb Z^\kappa),\mu\x\#,T_\phi)$ is rationally weakly mixing when $\kappa=1,2$.
\proclaim{Proposition 4.2\ \ ({\tt RWM} of special semiflows)}
\

Suppose that  the function $(h,\phi):\Om\to\Bbb R\x\Bbb Z^\kappa$ is  {\tt  aperiodic} in the sense that
$\nexists$ solution to
\begin{align*}&(h,\phi)=k+g-g\circ T,\ \ g:\Om\to\Bbb R\x\Bbb Z^\kappa\ \ \text{measurable},\\ & k:\Om\to \Bbb K\ \ 
\text{measurable, where $\Bbb K$ is a proper coset of $\Bbb R\x\Bbb Z^\kappa$}, 
\end{align*}
then for each $t>0$, $(X,\B,m,\Psi_t)$ is rationally weakly mixing.\endproclaim
 {\tt RWM}  entails a  kind of ``density local limit theorem'' and we begin the proof of proposition 4.2 with a one-sided version of ({\tt LLT}) for 
 the flow transformations.
\proclaim{Lemma 4.3: \ \ (Lower local limit)}
\

Suppose that  the function $(h,\phi):\Om\to\Bbb R\x\Bbb Z^\kappa$ is  {\tt  aperiodic} in the sense that
$\nexists$ solution to
\begin{align*}&(h,\phi)=k+g-g\circ T,\ \ g:\Om\to\Bbb R\x\Bbb Z^\kappa\ \ \text{measurable},\\ & k:\Om\to \Bbb K\ \ 
\text{measurable, where $\Bbb K$ is a proper coset of $\Bbb R\x\Bbb Z^\kappa$}, 
\end{align*}
then for $A\subset\Om$ a cylinder set and $I\subset [0,\min h]$ an interval,
\begin{align*}\tag{{\tt LLL}} \varliminf_{t\to\infty}t^{\frac{\kappa}2}\widehat{\Psi}_t(1_A\otimes 1_{\{0\}}\otimes 1_I)(\om,0,y)\ge 
\varkappa^{\frac{\kappa}2-1}f_X(0)\mu(A)|I|.
\end{align*}\endproclaim\demo{Proof  of ({\tt LLL})}

We have that
\begin{align*}\widehat{\Psi}_t(1_A\otimes 1_{\{0\}}\otimes 1_I)(\om,y,0)=\sum_{n=0}^\infty\T^n(f1_{[\phi_n=0,\ h_n\in I+t-y]})(\om)
\end{align*}
and so it suffices to prove
\begin{align*}\tag{\Letter} \lim_{M\to\infty}\lim_{t\to\infty}t^{\frac{\kappa}2}
\sum_{n=\frac{t}\varkappa \pm M\sqrt t}\T^n(&1_{A\cap [\phi_n=0,\ h_n\in I+t-y]})\\ &=\varkappa^{\frac{\kappa}2-1}f_X(0)\mu(A)|I|.
\end{align*}
\demo{Proof of (\Letter)}
\

 By the central limit theorem for $(h,\phi)$,
$$\frac1{\sqrt n}(\phi_n,h_n-\varkappa n)\xrightarrow[n\to\infty]{\text{\tiny distribution}}\ Z$$
where $Z=(X,Y)$ is non-singular normal with 
\bul $X$   centered, non-singular normal on $\Bbb R^\kappa$,\bul $Y$ centered normal on $\Bbb R$.

By the {\tt LLT} for $(h,\phi)$,
\begin{align*}
n^{\frac{\kappa+1}2}\T^n(1_A1_{[\phi_n=0,\ h_n\in I+n\varkappa+x_n\sqrt{n}]}) \ \underset{n\to\infty,\ |x_n|\le M}{\text{\Large $\approx$}}\ f_Z(0,x_n)\mu(A)|I|.
\end{align*}

Fix $t,\ M>0$, then
\begin{align*}&t=\varkappa n+x\sqrt{n}\ \ \text{with} \ \ x=x_{n,t}\in [-M,M]\ \ \iff\\ &
  n=\frac{t}\varkappa -\frac{x}{\varkappa^{\frac32}} \sqrt t\ \ \&\ \text{in this case}\\ &
 \T^n(f1_{[\phi_n=0,\ h_n\in I+t-y]})  \ \sim\ \frac{|I|}{n^{\frac{\kappa+1}2}}f_Z(0,x_{n,t})\\ &
\text{as}\ t,\ n\to\infty,\ |x_{n,t}|\le M.
\end{align*}

It follows that for fixed $M>0$ with $M':=\frac{M}{\varkappa^{\frac32}}$,
\begin{align*}
t^{\frac{\kappa}2}\sum_{n=\frac{t}\varkappa \pm M'\sqrt t}&\T^n(f1_{[\phi_n=0,\ h_n\in I+t-y]})\\ &
\ \underset{t\to\infty}{\text{\Large $\sim$}}\ \varkappa^{\frac{\kappa}2}n^{\frac{\kappa}2}\sum_{n=\frac{t}\varkappa \pm M'\sqrt t}\T^n(f1_{[\phi_n=0,\ h_n\in I+t-y]})
\\ & \ \underset{t\to\infty}{\text{\Large $\sim$}}\   \varkappa^{\frac{\kappa}2}\sum_{n=\frac{t}\varkappa \pm M'\sqrt t}\frac{|I|}{\sqrt{n}} f_Z(0,x_{n,t})
\end{align*}
Now,
$$\frac1{\sqrt{n}}-\frac1{\sqrt{n+1}}\sim\frac1{2n\sqrt{n}}$$ so
\begin{align*}x_{n,t}-x_{n+1,t}&=\frac{t-\varkappa n}{\sqrt{n}}-\frac{t-\varkappa (n+1)}{\sqrt{n+1})}
\\ &=\frac{\varkappa}{\sqrt{n+1}}+(t-\varkappa n)(\frac1{\sqrt{n}}-\frac1{\sqrt{n+1}})\\ &
=\frac{\varkappa}{\sqrt{n}}(1+ O(\frac1{\sqrt n})).
\end{align*}
Thus 

\begin{align*}\sqrt{t}\sum_{n=\frac{t}\varkappa \pm M'\sqrt t}&\T^n(f1_{[\phi_n=0,\ h_n\in I+t-y]})\\ &
\ \underset{t\to\infty}{\text{\Large $\approx$}}\  \varkappa^{\frac{\kappa}2}\sum_{n=\frac{t}\varkappa \pm M'\sqrt t}\frac{|I|}{\sqrt{n}} f_Z(0,x_{n,t})\\ &
\ \underset{t\to\infty}{\text{\Large $\approx$}}\  \varkappa^{\frac{\kappa}2-1} |I|\sum_{n=\frac{t}\varkappa \pm M'\sqrt t}(x_{n+1,t}-x_{n,t}) f_Z(0,x_{n,t})\\ &
\xrightarrow[t\to\infty]{}\  \varkappa^{\frac{\kappa}2-1} |I|\int_{[-M',M']}f_Z(0,x)dx \\ &\xrightarrow[M'\to\infty]{}  \varkappa^{\frac{\kappa}2-1} |I|f_X(0).\ \ \CheckedBox\end{align*}
\demo{Proof of proposition 4.2} This follows from pointwise dual ergodicity and (LLL) via proposition 3.3 of \cite{RatWM}.\ \Checkedbox
\subsection*{Remark}\ \ A stronger version of  lemma 4.3 would be  the  {\tt local limit theorem}:
\sms {\it For $A\subset\Om$ a cylinder set and $I\subset [0,\min h]$ an interval,}
\begin{align*}\tag{{\tt LLT}} \lim_{t\to\infty}t^{\frac{\kappa}2}\widehat{\Psi}_t(1_A\otimes 1_{\{0\}}\otimes 1_I)(\om,0,y)=
\varkappa^{\frac{\kappa}2-1}f_X(0)\mu(A)|I|.
\end{align*}
 It is not hard to show that  special semiflows satisfying ({\tt LLT}) enjoy the stronger property of {\tt Krickeberg mixing} (as in \cite{Krickeberg}).
Note that here, in the notation of lemma 4.3, ({\tt LLT}) is equivalent to
\begin{align*}\tag{\Large\bell} \varlimsup_{M\to\infty}\varlimsup_{t\to\infty}t^{\frac{\kappa}2}
\sum_{n\ge 1,\ |n-\frac{t}\varkappa|\ge M\sqrt t}\T^n(1_{[\phi_n=0]})=0.
\end{align*}
We do not know whether this necessarily holds under the assumptions of lemma 4.3, but we'll see  
in the next section that the geodesic flows of Abelian covers of compact hyperbolic surfaces have this $\varlimsup$ finite when we show (\dsrailways) (on page \pageref{choochoo}).
\

For work on ({\tt LLT}) for flows, see \cite{Iwata} $\&$ \cite{Waddington}.
\section*{ \S5 hyperbolic geodesic flows}
\subsection*{Definitions}
The  {\it hyperbolic plane} is
$\Bbb H:=\{z=u+iv\in\Bbb C:|z|<1\}$ equipped with the arclength element
$ds(u,v):={2\sqrt{du^2+dv^2}\over 1-u^2-v^2}$ and the area element
$dA(u,v):={4dudv\over (1-u^2-v^2)^2}.$

The {\it hyperbolic distance} between $x,y\in \Bbb H$ is 
$$\rho(x,y):=\ \inf\,\{\int_{\gamma}ds:\ \gamma\ \text{is an arc joining}\
 x\ \text{and}\ y\}
= 2\tanh^{-1}{|x-y|\over |1-\overline xy|}.$$
This $\inf$ is achieved by an arc of a {\it geodesic} in $(\Bbb H,\rho)$. Note that $\Bbb H$ is also known as  the {\tt disk model} for the hyperbolic plane. 
These geodesics are diameters of $\Bbb H$ and circles orthogonal
 to $\partial \Bbb H$.\footnote{ It is sometimes convenient to consider the conformally equivalent {\tt upper half plane model} $\widehat{\Bbb H}:=\{x+iy:\ x,y\in\Bbb R,\ y>0\}$
 where the geodesics are vertical lines and circles orthogonal
 to $\Bbb R$.}

\par The isometries $\isom(\Bbb H,\rho)$ of $(\Bbb H,\rho)$ are  the  M\"obius transformations
 and their complex conjugates.
\par If $g$ is an isometry of $\Bbb H$, then $A\circ g\equiv A$.

\par The space of {\it line elements} of $\Bbb H$ is $UT(\Bbb H)=\Bbb H\x\Bbb T$.

\par The {\it geodesic flow transformations} $\v^t$ are defined on
$\Bbb H\x\Bbb T$ as follows.
To each line element $\om$
there corresponds a unique directed geodesic passing through $x(\om)$ whose
directed tangent at $x(\om)$
makes an angle $\theta(\om)$ (with the radius $(0,1)$).
\par If $t>0$, the point $x(\v^t\om)$ is the unique point on
the geodesic at distance
 $t$ from $x(\om)$ in the direction of the geodesic, and if
$t<0$, the point $x(\v^t\om)$ is the unique point on the geodesic at distance $-
t$ against the direction of the geodesic.
\par The angle $\theta(\v^t\om)$ is the angle made by the directed
tangent to the geodesic at the point $x(\v^t\om)$.
\smallskip
\par There is an important involution $\chi :\Bbb H\x\Bbb T\to \Bbb H\x\Bbb T$,  of
direction reversal: $x(\chi \om)=x(\om)$ and $\th(\chi\om)=\th(\om)+\pi$. Here $\om=(x,\th)=(x(\om),\th(\om))$.
\par The isometries act on $\Bbb H\x\Bbb T$ (as differentiable maps) by
$$g(\om)=(g(x(\om)),\theta(\om)+\arg g'(x(\om))$$
and it is not hard to see that $\chi  g=g\chi $ and $\v^tg=g\v^t$.
\par Both the geodesic flow, the involution and the
isometries preserve the measure
\par $dm(x,\theta)=dA(x)d\theta$ on $\Bbb H\x\Bbb T.$

\

\par Let $\G$ be a discrete subgroup of  $\text{\tt Isom}\,(\Bbb H)$ (aka {\it Fuchsian group}), then $M=\Bbb H/\G$ is an hyperbolic surface and any hyperbolic
surface is isometric to one of this form.
\

The space of line elements of $M=\Bbb H/\G$ is
$UT(M):=M\x\Bbb T =(\Bbb H\x\Bbb T )/\G$
and the geodesic flow transformations on $UT(M)$ are defined by
$$\v^{M}_t\G(\om):=\G \v^t(\om).$$

\par Let $\pi_\G:\Bbb H\to M,\ \overline\pi_\G:\Bbb H\x\Bbb T\to UT(M)$ be the projections
 $\pi_\G(z)=\G z,\ \overline\pi_\G(\om)=\G\om$,
and let $F$ be a {\it fundamental domain}
for $\G$ in $\Bbb H$, e.g.
$$F^o:=\{x\in \Bbb H:\rho(y,x)<\rho(\g(y),x)\ \forall\
\g\in\G\setminus\{e\}\},\ \ y\in \Bbb H,$$
then $\pi_\G$ and $\overline\pi_\G$ are 1-1 on $F$ and $F\x\Bbb T$,
and so the measures $A_{|F}$ and $m_{|F}$ induce measures $A_\G$ and $m_\G$ on
$M=\Bbb H/\G$ and $UT(M)=\Bbb H/\G\x\Bbb T$ respectively.
\subsection*{Basics}
\par It is known that  for $M=\Bbb H/\G$,
\sbul $\v^M$ is either totally dissipative, or
conservative and ergodic  (E. Hopf \cite{Hopfbook}), 

\begin{align*}\tag*{$\bullet$}\ \v^M\ \ \text{is conservative iff}\ \ a_\G(t):=
\sum_{\g\in\G,\ \rho(x,\g(x))\le t}e^{-\rho(x,\g(x))}\xrightarrow[t\to\infty]{}\ \infty
\end{align*}
(E. Hopf \cite{Hopf} $\&$ M. Tsuji \cite{Tsuji}).

\

Moreover, any conservative
 $\v^M$ is:
 \bul rationally ergodic  with  return sequence
 $\propto\  a_\G(t)$ (\cite{A-S}, see also \cite{IET} chapter 7);
 \bul weakly mixing (\cite{Roblin}).
  Note that a flow is weakly mixing iff all its transformations are ergodic.
  \
  
   All transformations of a rationally ergodic, weakly mixing flow are necessarily rationally ergodic.

\subsection*{Abelian covers of compact surfaces}
\

\par Let $M=\Bbb H/\G$ be a compact, hyperbolic surface, let $\v^M:UT(M)\to UT(M)$ denote the
geodesic flow and let $\chi:UT(M)\to UT(M)$ be the involution of
direction reversal.

\par Now let $\kappa\ge 1\ \&$ let
$V=V^{(\kappa)}$ be a $\Bbb Z^\kappa$-cover of $M$ that is $V$ is a complete
hyperbolic surface equipped with a covering map $p:V\to M$ so that
$\exists$ a monomorphism $\g:\Bbb Z^\kappa\to\text{\tt Isom}\,(V^{(\kappa)})$, such that
for $y\in V,\ p^{-1}\{p(y)\}=\{\g(n)y:\ n\in\Bbb Z^\kappa\}$. 
\

Rees showed in \cite{Rees} that $\v^{V^{(\kappa)}}$ is conservative when $\kappa=1,2$ and dissipative when $\kappa\ge 3$.
\

In this section we prove
\proclaim{Theorem 5.1}\ \ The  geodesic flow transformations ${\v_t^{V^{(\kappa)}}}\ \ (t>0)$ are
\bul rationally weakly mixing when $\kappa=1, 2$;
\bul  $2$-recurrent when $\kappa=1$ and $2$-dissipative when $\kappa=2$.\endproclaim
\demo{Proof of rational weak mixing}
 
\par Let $M=\Bbb H/\G$ be a compact, hyperbolic surface (with $\G$ is the corresponding, cocompact, Fuchsian group) and let $\v_M:UT(M)\to UT(M)$ denote the
geodesic flow on $UT(M)$ (the unit tangent bundle) and let $\chi:UT(M)\to UT(M)$ be the involution of
direction reversal.
\par As $\v^M$ is an Anosov flow, by Bowen's theorem (\cite{Bowen}), there is a special flow
$\Phi:\Bbb R\to\text{\tt PPT}(Y,\mathcal C,\nu)$ and $\pi:Y\to UT(M)$  a continuous, measure theoretic isomorphism satisfying $\v^M\circ\pi=\pi\circ\Phi$.
Here $\text{\tt PPT}(Y,\mathcal C,\nu)$ denotes the collection of invertible probability preserving transformations of the probability space
$(Y,\mathcal C,\nu)$ equipped with the weak operator topology.
\

Here:
\begin{align*}
 &Y=\{(x,t)\in\Om\x\Bbb R_+:\ 0\le x<h(x)\}\ \ \ \text{where}\\ &
 \Om\subset S^\Bbb Z\ \ \text{a transitive {\tt SFT},\ $h:\Om\to\Bbb R_+$  H\"older};\\ & \mathcal C=\B(Y),\ \ \ \nu(A\x B)=c^{-1}\mu(A)\text{\tt Leb}\,(B)
 \\ &\text{$\mu\in\mathcal P(\Om)$ Gibbs as in \cite{BowenBook};\ $c:=\int_\Om hd\mu$}\
 \& \ \Phi_t(x,y)=(T^nx,y+t-h_n(x))\\ &\  \text{where $T$ is the shift and $n=n_t(x,y)$ is so that}\\ & h_n(x):=\sum_{k=0}^{n-1}h(T^kx)\le y+t<h_{n+1}(x).
\end{align*}

\

\par By Rees' refinement, $(\Om,T,\mu),\ h\ \&\ \pi$ can be chosen so that
\bul  $S$ is a finite, symmetric generator set of $\G$ and the elements of $\Om$ code the geodesics in $M=\Bbb H/\G$,
\sbul   $(\Om,T,\mu)$ is topologically mixing, \sbul  $h(\dots,x_{-1},x_0,x_1,\dots)=h(x_1,x_2,\dots)$ and \sbul  $\chi(\pi\Sigma)=\pi\Sigma$.

\par Now let $\kappa\ge 1\ \&$ let
$V$ be a $\Bbb Z^\kappa$-cover of $M$ that is $V$ is a complete
hyperbolic surface equipped with a covering map $p:V\to M$ so that
there exists a monomorphism $\g:\Bbb Z^\kappa\to\text{\tt Isom}\,(V)$, such that
for $y\in V,\ p^{-1}\{p(y)\}=\{\g(n)y:\ n\in\Bbb Z^\kappa\}$. 
\

The corresponding tangent map, also denoted  $p:UT(V)\to UT(M)$ is equivariant with the geodesic flows and their direction reversal involutions.

We have that $V^{(\kappa)}\cong\Bbb H/\G_0$ where  the corresponding  Fuchsian group $\G_0=\text{\tt Ker}\,\Theta$ for $\Theta:\G\to\Bbb Z^\kappa$ a 
surjective homomorphism.
\

The corresponding $\Bbb Z^\kappa$-extension of $\Phi:\Bbb R\to\text{\tt PPT}(Y,\mathcal C,\nu)$  is the
special flow $\Psi:\Bbb R\to\text{\tt MPT}\xbm$  where
\begin{align*}
 &X=\{(x,n,t)\in\Om\x\Bbb Z\x\Bbb R_+:\ 0\le x<h(x)\},\\ & 
 \ m(A\x B\x C)=\mu(A)\text{\tt Leb}\,(B)\#(C)
\ \ \text{where $\#$ is counting measure,  }
\\ &
   \Psi_t(x,y,z)=(\Phi_t(x,y),z+\phi_n(x)) \ \ \text{where}\\ & \Phi_t(x,y)=(T^nx,y+t-h_n(x)) \ \ \&\ \phi(\om_1,\om_2,\dots)=\Theta(\om_1).
\end{align*}
There is a continuous, measure theoretic isomorphism $\Pi:X\to TV$  satisfying 
$$p\circ\Pi\equiv\pi,\ \ \v_V\circ\pi=\pi\circ\Psi\ \ \&\ \ \Pi(x,t,n):=\g(n)\Pi(x,t,0).$$

\ As in  \cite{Sharp}, $V$ is {\it homologically full} in the sense that as $t\to\infty,\ \exists$ exponentially many closed geodesics of length $\le t$ in each homology class.
\

\sms Therefore,  by the lemma in \cite{Solomyak}  the function $(h,\phi):\Om\to\Bbb R\x\Bbb Z^\kappa$ is  {\tt  aperiodic} in the sense that
$\nexists$ solution to
\begin{align*}&(h,\phi)=k+g-g\circ T,\ \ g:\Om\to\Bbb R\x\Bbb Z^\kappa\ \ \text{measurable},\\ & k:\Om\to \Bbb K\ \ \text{measurable, where $\Bbb K$ is a proper coset of $\Bbb R\x\Bbb Z^\kappa$}.
\end{align*}
Thus $\Psi:\Bbb R\to\text{\tt MPT}\xbm$  is a two-sided version of a special semiflow satisfying the assumptions of proposition 4.2
and rational weak mixing follows.\ \ \ \Checkedbox
\demo{Proof of $2$-recurrence $\&$ $2$-dissipation}

For $x\in \Bbb H$, and $\e>0$, set
$$N_\rho(x,\e)=\{y\in \Bbb H:\rho(x,y)<\e\},\ \ \D (x,\e):=N_\rho(x,\e)\x\Bbb T.$$
To prove the theorem, we show first that
\

sets of form $\D=\D(x,\e)$ are admissible and satisfy 
\begin{align*}
 \tag{\dsrailways}\label{choochoo}m(\D\cap\v_{V^{(d)}}^{-t}\D)\ \asymp\ \frac1{t^{\frac{d}2}}.
\end{align*}
The lower bound in (\dsrailways) follows from lemma 4.3 applied to the flow $(X,\B,m,\Psi)$ as above. 
\

We now proceed to establish the upper bound. The following analytic geometry lemmas are valid for any Fuchsian group $\G$ with 
$X_\G=UT(\Bbb H/\G)=\Bbb H/\G\x\Bbb T$.

\

\proclaim{ Analytic geometry lemma I}: 
 \
 
 For $x\in X_\G\ \ \&\ \e>0$ small enough:
\begin{align*}&\tag{i}m(\D(x,\e)\cap \phi_\G^{-t}\D(x,\e))\ \ll\ \sum_{\g\in\G,\ \rho(x,\g(x)=t\pm 2\e}e^{-\rho(x,\g(x))}\\ &\tag{ii}
m(\D(x,\e)\cap \phi_\G^{-t}\D(x,\e))\gg \sum_{\g\in\G,\ \rho(x,\g(x)=t\pm \frac{\e}2}e^{-\rho(x,\g(x))}.\end{align*}
\endproclaim
\demo{ Proof}\ \ We have that
\begin{align*} m(\D(x,\e)\cap\v_\G^{-s}\D(x,\e)) &=
\int_{\D(x,\e)}1_{\D(x,\e)}\circ \v_\G^{s}dm_\G\\ &=
\int_{N_\rho(x,\e)}\Phi(s;z)dA(z)\end{align*}
where
$$\Phi(s;z):=
\sum_{\g\in\G}\int_{\Bbb T}1_{\g N_\rho(x,\e)\x\Bbb T}\circ \v_\G^s
(z,\theta)d\theta.$$
Set $\v_z(\om)={z+\om\over 1+\overline z\om}$. Using
$\v_z\v^t=\v^t\v_z$, and $\v^t(0,\th)=(\tanh {t\over 2} e^{2\pi i\th},\th)$,
we have
\begin{align*}  \Phi(s;z) & =
\sum_{\g\in\G}\int_{\Bbb T}1_{\g N_\rho(x,\e)\x\Bbb T}\circ \v_\G^s(z,\theta)d\theta\\
 &=
 \sum_{\g\in\G}\int_{\Bbb T}
1_{\v_z^{-1}\g N_\rho(x,\e)\x\Bbb T}\circ \v_\G^s(0,\theta)d\theta\\ & =
 \sum_{\g\in\G}\int_{\Bbb T}1_{\v_z^{-1}\g N_\rho(x,\e)}
(\tanh(\tfrac{s}{2}) e^{2\pi i\theta})d\theta\\ &=
\sum_{\g\in\G}\int_{\Bbb T}1_{ N_\rho(\v_z^{-1}\g(x),\e)}
(\tanh(\tfrac{s}{2}) e^{2\pi i\theta})d\theta\\ &=
\sum_{\g\in\G}|J(\v_z^{-1}\g(x),\e)|
\end{align*}
where $J(w,\eta)\subset\Bbb T$ is the interval
$$J(w,\eta):=\{\th\in\Bbb T:\ \tanh(\tfrac{s}{2}) e^{2\pi i\theta}\in N_\rho(w,\eta)\}$$ and
$|J(w,\eta)|$ is its length.

We have that $|J(w,\eta)|>0$ iff $\rho(0,w)=s\pm\eta$. 
\

Thus
$|J(\v_z^{-1}\g(x),\e)|>0$ iff
$$\rho(z,\g(x))=\rho(0,\v_z^{-1}\g(x))=s\pm\e$$
and
\begin{align*}\Phi(s;z)= \sum_{\g\in\G,\ \rho(z,\g(x))=s\pm\e}|J(\v_z^{-1}\g(x),\e)|. 
\end{align*}

\par Next, for $w\in\Bbb H\ \&\ \eta>0$, we consider the angle interval
subtended by $N_\rho(w,\eta)$ at $0\notin N_\rho(w,\eta)$,
$$\L(w,\eta):=\{\theta\in [0,2\pi]:\exists\ r\in (0,1)\ \ni\
\rho(w,re^{i\theta})<\eta\}.$$
We note that
$$\L(w,\eta)=\{\theta\in [0,2\pi]:\|\theta-\arg w\|<\sin^{-1}\(
{(1-|w|^2)\tanh{\eta\over 2}\over |w|(1-\tanh^2{\eta\over 2})}\)\},$$
where $\|\theta\|:=\theta\wedge(2\pi-\theta),\ \ \theta\in [0,2\pi)$.
 This is because
$$N_\rho(w,\eta)=B\({(1-\d^2)w\over 1-\d^2|w|^2},{\d(1-|w|^2)\over
1-\d^2|w|^2}\)$$
where $B(x,r)$ is the Euclidean ball of radius $r$ and $\d=\tanh{\eta\over 2}$. 
 Thus as  $|w|\to 1\ \& \ \eta\to 0$,
\begin{align*} |\L(w,\eta)|&=2\sin^{-1}\({(1-|w|^2)\tanh{\eta\over 2}\over |w|(1-
\tanh^2{\eta\over 2})}\)\\
& \sim
\eta(1-|w|^2)\\& \sim\  \eta e^{-\rho(0,w)}
\end{align*}
Thus,
\begin{align*}\Phi(s;z)&= \sum_{\g\in\G,\ \rho(z,\g(x))=s\pm\e}|J(\v_z^{-1}\g(x),\e)|\\ &\le
\sum_{\g\in\G,\ \rho(z,\g(x))=s\pm\e}|\L(\v_z^{-1}\g(x),\e)|\\ &\ll
\sum_{\g\in\G,\ \rho(z,\g(x))=s\pm\e}e^{-\rho(z,\g(x))}\\ &\le \sum_{\g\in\G,\ \rho(x,\g(x))=s\pm 3\e}e^{-\rho(x,\g(x))}.
\end{align*}
whence
\begin{align*} m(\D(x,\e)\cap\v_\G^{-s}\D(x,\e)) &=
\int_{N_\rho(x,\e)}\Phi(s;z)dA(z)\\ &\ll\sum_{\g\in\G,\ \rho(x,\g(x))=s\pm 2\e}e^{-\rho(x,\g(x))}.\ \ \CheckedBox{\rm(i)} \end{align*}
Next, to establish (ii) note that $\exists\ \zeta>0$ so that
$$\rho(0,w)=s\pm\frac{3\eta}4\ \ \Lra\ \ |J(w,\eta)|>\zeta|\L(w,\eta)|.$$
Thus, for $z\in N_\rho(x,\frac{\e}4)$
\begin{align*}\Phi(s;z)&\ge  \sum_{\g\in\G,\ \rho(z,\g(x))=s\pm\frac{3\e}4}|J(\v_z^{-1}\g(x),\e)|\\ &\ge
\zeta\sum_{\g\in\G,\ \rho(z,\g(x))=s\pm\frac{3\e}4}|\L(\v_z^{-1}\g(x),\e)|\\ &\gg
\sum_{\g\in\G,\ \rho(z,\g(x))=s\pm\frac{3\e}4}e^{-\rho(z,\g(x))}\\ &\ge \sum_{\g\in\G,\ \rho(x,\g(x))=s\pm \frac{\e}2}e^{-\rho(x,\g(x))}
\end{align*}
whence
\begin{align*} m(\D(x,\e)\cap\v_\G^{-s}\D(x,\e)) &\ge
\int_{N_\rho(x,\frac{\e}4)}\Phi(s;z)dA(z)\\ &\gg\sum_{\g\in\G,\ \rho(x,\g(x))=s\pm \frac{\e}2}e^{-\rho(x,\g(x))}.\ \ \CheckedBox{\rm(ii)} \end{align*}

\proclaim{ Analytic geometry lemma II}: 
 \
  
 For $x\in X_\G\ \ \&\ \e>0$ small enough:

\begin{align*}
m(\bigcap_{j=0}^p\phi_\G^{-\sum_{i=0}^js_i}&\D(x,\e))\le M_p\prod_{k=1}^pm(\D(x,4\e)\cap \phi_\G^{-s_k}\D(x,4\e))\\& \forall\ s_0=0,\ s_1,\dots,s_p>0. 
\end{align*}\endproclaim
\demo{ Proof}

\

For $\u \g\in \G^p$ (resp. $\u t\in \Bbb R^p$) we denote its coordinates by
$\g_k$ (resp. $t_k$), $k=1,...,p$. Let
$$ I_p=\{\u t\in \Bbb R^p: 0<t_1<...<t_p\}.$$

Let $\e>0$ be fixed and $\D=N\x \Bbb T$ as
before, where $N=N_{\rho}(x,\e)$. We assume $\e$ to be sufficiently small.

First observe that
\begin{align*}u(p,\u t)&:=m(\bigcap_{j=0}^p\phi_\G^{-t_{j}}\D)\\
&=\int_{\D}\prod_{j=1}^p
1_{\D}\circ\v_\G^{t_j} dm_{\G}\\
& =\sum_{\u \g\in\G^p}\int_{\D}
\prod_{j=1}^p1_{\g_j\D}\circ\v^{t_j}dm\\
& =\sum_{\u \g\in\G^p}\int_{N}\int_0^{2\pi}
\prod_{j=1}^p
1_{\g_j N\x\Bbb T}\circ\v^{t_j}(z,\theta)d\theta dA(z)\\ &=\int_N \psi_p(\u t,z) A(dz)
\end{align*}
where
$$\psi_p(\u t,z):=\sum_{\u \g\in\G^p}\int_0^{2\pi}
\prod_{j=1}^p1_{\g_j N\x\Bbb T}\circ\v^{t_j}(z,\theta) d\theta.$$

Next,
\begin{align*} \psi_p(t,z)=\sum_{\u \g\in\G^p}\int_0^{2\pi}
\prod_{j=1}^p1_{\v_z^{-1}\g_j N}(\tanh t_j e^{i\theta})
d\theta
 \end{align*}

For $\u t\in I_p$, , let $t_0=0$, let
$$s_{k+1}=t_{k+1}-t_k\ \ (0\le k\le p-1)$$  and let
$$\G_0(\u t):= \{\u \g\in\G^p:
\int_0^{2\pi}
\prod_{j=1}^p1_{\v_z^{-1}\g_j N}(\tanh t_j e^{i\theta})
d\theta>0\}.$$
If  $\u\g=(\g_1,\dots,\g_p)\in\G_0(\u t)$, then 
$$\exists\ \th\in [0,2\pi)\ \text{with}\ \ \rho(\om(\v_{t_k}(z,\th)),\g_k(x))<\e\ \forall\ 1\le k\le p$$ whence
\begin{align*}\rho(\g_k(x),\g_{k+1}(x))&=\rho(\om(\v_{t_k}(z,\th)),\om(\v_{t_{k+1}}(z,\th)))\pm2\e\\ &=s_k\pm2\e.\end{align*}

Setting $w_0=z=\g_0(x)$,
\begin{align*}
 \rho(z,\g_p(x))&=t_p\pm\e\\ &=\sum_{k=0}^{p-1}s_k\pm\e\\ &=\sum_{k=0}^{p-1}\rho(\g_{k+1}(x),\g_k(x))\pm (2p+1)\e.
\end{align*}

Thus
\begin{align*}\psi_p(\u t,z)&=\sum_{\u \g\in\G^p(\u t)}\int_0^{2\pi}
\prod_{j=1}^p1_{\v_z^{-1}\g_j N}(\tanh t_j e^{i\theta})
d\theta\\ &\le \sum_{\u \g\in\G^p(\u t)}\int_0^{2\pi}
1_{\v_z^{-1}\g_p N}(\tanh t_j e^{i\theta})
d\theta\\ &\le \sum_{\u \g\in\G^p(\u t)}|\L(\v_z^{-1}\g_p(x),\e)|\ \ \text{\scriptsize $\because\ \tanh t_j e^{i\theta}\in \v_z^{-1}\g_j N\ \Rightarrow\ \th\in\L(\v_z^{-1}\g_p(x),\e)$},
\\ &\ll \sum_{\u \g\in\G^p(\u t)}e^{-\rho(0,\v_z^{-1}\g_px)}\\ &=
\sum_{\u \g\in\G^p(\u t)}e^{-\rho(z,\g_px)}\\ &\ll \sum_{\u \g\in\G^p(\u t)}\prod_{k=0}^{p-1}e^{-\rho(\g_k(x),\g_{k+1}(x))}\\ &\le
\prod_{k=0}^{p-1}\sum_{ \g\in\G,\ \rho(z,\g(x))=s_k\pm 2\e}e^{-\rho(x,\g(x))}\\ &\ll
\prod_{k=0}^{p-1}m(\D(x,4\e)\cap\v_\G^{-s_k}\D(x,4\e)\end{align*}
and
\begin{align*}
m(\bigcap_{j=0}^p\phi_\G^{-t_j}\D(x,\e))&=\int_N\psi_p(\u t,z)dA(z)\\ &\ll
\prod_{k=1}^pm(\D(x,4\e)\cap \phi_\G^{-s_k}\D(x,4\e)).\ \ \ \ \CheckedBox 
\end{align*}
To complete the proof, we use  a {\tt word metric observation} from \cite{Rees}.
\subsection*{Word length}
\

Define the $S$-{\it word length} of $\g\in\G$ by
$$\ell(\g)=\ell_S(\g):=\min\,\{N\ge 1:\ \exists\ c_1,c_2,\dots,c_N\in S,\ \g= c_1c_2\dots c_N\}.$$
This gives rise to the {\it word metric} $d_\ell$ on $\G$ given by 
$$d_\ell(\b,\g):=\ell(\g b^{-1}).$$

As before, a set of form 
$$C=[c_1,c_2,\dots,c_n]:=\{x\in \Om:x_k=c_k\ \forall\ 1\le k\le n\}$$
is called a {\it cylinder of length $n$}. Let
$$\mathcal C_n:=\{\text{cylinders of length $n$}\}.$$
To each $C=[c_1,c_2,\dots,c_n]\in\mathcal C_n$ corresponds $\g=\g_C:=c_1c_2\dots c_n\in\G$ with $\ell(\g_C)=n$.

It is shown  in \cite{Rees} that $\exists\ M=M_\G>0$ so that 
\begin{align*}\tag{\Wheelchair}& \text{ (i)}\ \ \rho(\g(0),\b(0))\ =\ \ M^{\pm}d_\ell(\g,\b),\\ &
\text{ (ii)}\ \ \mu(C)=M^{\pm 1}e^{-\rho(0,\g(0))}\\ &
\text{where}\ \emptyset\ne C=[c_1,c_2,\dots,c_n]\subset\Om\ \&\ \g=\g_C:=c_nc_{n-1}\dots c_2c_1. 
\end{align*}

Thus,
\begin{align*}
 \sum_{\g\in\G_0,\ \ell(\g)=n}e^{-\rho(0,\g(0))}
 & =M^{\pm}\sum_{\g=c_1c_2\dots c_n\in\G_0,\ \ell(\g)=n}\mu([c_1,\dots,c_n])\ \ \text{by (\Wheelchair)(ii)}\\ &=
 \sum_{C\in\mathcal C_n,\ \Theta(\g_C)=0}\mu([c_1,\dots,c_n])\\ &=\mu([\phi_n=0])\\ &\asymp \frac1{n^{\frac{\kappa}2}}\ \ \ \text{by ({\tt LLT})}.
\end{align*}
Fix $t,\ K>0$. Suppose that  $\frak g\in\G_0\ \&\ \rho(0,\frak g(0))=t\pm K$. Let $\ell(\frak g)=N$.
If $\g\in\G,\ \rho(0,\g(0))=t\pm K$, then $\rho(\frak g(0),\g(0))<2K$, whence $d_\ell(\frak g,\g)<2MK$ and
$$\ell(\g)=N\pm 2MK.$$
Using this and (\Wheelchair)(i)
\begin{align*}\sum_{\g\in\G_0,\ \rho(0,\g(0))=t\pm K}e^{-\rho(0,\g(0))}&
\le \sum_{n=N\pm 2MK}\sum_{C\in\mathcal C_n,\ \Theta(\g_C)=0}e^{-\rho(0,\g(0))}\\ &\le 
 M\sum_{n=N\pm 2MK}\sum_{C\in\mathcal C_n,\ \Theta(\g_C)=0}\mu([\phi_n=0])\\ &\ll \frac1{N^{\frac{\kappa}2}}\asymp\frac1{t^{\frac{\kappa}2}}.
\end{align*}
This is the upper estimation in (\dsrailways) and proves admissibility of  sets of form $\D(x,\e)$.
It follows that $\v_{V^{(2)}}$ is $1$-nice and $2$-dissipative. 
\

It follows from ({\tt LLL}) that the representing semiflow of $\v_{V^{(1)}}$ is $2$-nice, whence also $\v_{V^{(1)}}$. By theorem 2.2, $\v_{V^{(1)}}$ is 
subsequence $2$-rationally ergodic, whence $2$-recurrent.\  \ \Checkedbox 
\subsection*{Higher dimensional theorem 5.1}\ \ In the interest of simplicity, we stated and proved theorem 5.1 for surfaces.
In fact the analogous statements hold for  geodesic flows of hyperbolic manifolds of arbitrary dimension with constant negative curvature.
The proof is the same. The geodesic flow of a compact manifold with constant negative curvature is an Anosov flow. The  results of \cite{Bowen}
\cite{Rees}, \cite{Sharp} and \cite{Solomyak} apply   and the analytic geometry lemmas have analogous multidimensional versions.


\begin{thebibliography}{10}

\bibitem{RatErg}
Jon Aaronson.
\newblock Rational ergodicity and a metric invariant for {M}arkov shifts.
\newblock {\em Israel J. Math.}, 27(2):93--123, 1977.

\bibitem{IET}
Jon Aaronson.
\newblock {\em An introduction to infinite ergodic theory}, volume~50 of {\em
  Mathematical Surveys and Monographs}.
\newblock American Mathematical Society, Providence, RI, 1997.

\bibitem{RatWM}
Jon Aaronson.
\newblock Rational weak mixing in infinite measure spaces.
\newblock {\em Ergodic Theory and Dynamical Systems}, 33:1611--1643, 12 2013.

\bibitem{mode}
Jon Aaronson and Manfred Denker.
\newblock The {P}oincar\'e series of {$\Bbb C\setminus\Bbb Z$}.
\newblock {\em Ergodic Theory Dynam. Systems}, 19(1):1--20, 1999.

\bibitem{AD}
Jon Aaronson and Manfred Denker.
\newblock Local limit theorems for partial sums of stationary sequences
  generated by {G}ibbs-{M}arkov maps.
\newblock {\em Stoch. Dyn.}, 1(2):193--237, 2001.

\bibitem{ADSZ}
Jon. Aaronson, Manfred Denker, Omri Sarig, and Roland Zweim{\"u}ller.
\newblock Aperiodicity of cocycles and conditional local limit theorems.
\newblock {\em Stoch. Dyn.}, 4(1):31--62, 2004.

\bibitem{A-N}
Jon Aaronson and Hitoshi Nakada.
\newblock Multiple recurrence of {M}arkov shifts and other infinite measure
  preserving transformations.
\newblock {\em Israel J. Math.}, 117:285--310, 2000.

\bibitem{A-S}
Jon Aaronson and Dennis Sullivan.
\newblock Rational ergodicity of geodesic flows.
\newblock {\em Ergodic Theory Dynam. Systems}, 4(2):165--178, 1984.

\bibitem{Bowen}
Rufus Bowen.
\newblock Symbolic dynamics for hyperbolic flows.
\newblock {\em Amer. J. Math.}, 95:429--460, 1973.

\bibitem{BowenBook}
Rufus Bowen.
\newblock {\em Equilibrium states and the ergodic theory of {A}nosov
  diffeomorphisms}, volume 470 of {\em Lecture Notes in Mathematics}.
\newblock Springer-Verlag, Berlin, revised edition, 2008.
\newblock With a preface by David Ruelle, Edited by Jean-Ren{\'e} Chazottes.

\bibitem{G-L}
Adriano Garsia and John Lamperti.
\newblock A discrete renewal theorem with infinite mean.
\newblock {\em Comment. Math. Helv.}, 37:221--234, 1962/1963.

\bibitem{G-H}
Y.~Guivarc'h and J.~Hardy.
\newblock Th\'eor\`emes limites pour une classe de cha\^\i nes de {M}arkov et
  applications aux diff\'eomorphismes d'{A}nosov.
\newblock {\em Ann. Inst. H. Poincar\'e Probab. Statist.}, 24(1):73--98, 1988.

\bibitem{HH}
Hubert Hennion and Lo{\"{\i}}c Herv{\'e}.
\newblock {\em Limit theorems for {M}arkov chains and stochastic properties of
  dynamical systems by quasi-compactness}, volume 1766 of {\em Lecture Notes in
  Mathematics}.
\newblock Springer-Verlag, Berlin, 2001.

\bibitem{Hopfbook}
E.~Hopf.
\newblock {\em Ergodentheorie}.
\newblock Number v. 5, no. 2 in Ergebnisse der Mathematik und ihrer
  Grenzgebiete, 5. Bd. Julius Springer, 1937.

\bibitem{Hopf}
Eberhard Hopf.
\newblock Ergodic theory and the geodesic flow on surfaces of constant negative
  curvature.
\newblock {\em Bull. Amer. Math. Soc.}, 77:863--877, 1971.

\bibitem{Iwata}
Yukiko Iwata.
\newblock A generalized local limit theorem for mixing semi-flows.
\newblock {\em Hokkaido Math. J.}, 37(1):215--240, 02 2008.

\bibitem{Krickeberg}
Klaus Krickeberg.
\newblock Strong mixing properties of markov chains with infinite invariant
  measure.
\newblock In {\em Proceedings of the Fifth Berkeley Symposium on Mathematical
  Statistics and Probability, Volume 2: Contributions to Probability Theory,
  Part 2}, pages 431--446, Berkeley, Calif., 1967. University of California
  Press.

\bibitem{PPS}
F.~{Paulin}, M.~{Pollicott}, and B.~{Schapira}.
\newblock {Equilibrium states in negative curvature}.
\newblock {\em ArXiv e-prints}, November 2012.

\bibitem{Rees}
Mary Rees.
\newblock Checking ergodicity of some geodesic flows with infinite {G}ibbs
  measure.
\newblock {\em Ergodic Theory Dynamical Systems}, 1(1):107--133, 1981.

\bibitem{Roblin}
Thomas Roblin.
\newblock Sur l'ergodicit\'e rationnelle et les propri\'et\'es ergodiques du
  flot g\'eod\'esique dans les vari\'et\'es hyperboliques.
\newblock {\em Ergodic Theory Dynam. Systems}, 20(6):1785--1819, 2000.

\bibitem{Sharp}
Richard Sharp.
\newblock Closed orbits in homology classes for {A}nosov flows.
\newblock {\em Ergodic Theory Dynam. Systems}, 13(2):387--408, 1993.

\bibitem{Solomyak}
Rita Solomyak.
\newblock A short proof of ergodicity of {B}abillot-{L}edrappier measures.
\newblock {\em Proc. Amer. Math. Soc.}, 129(12):3589--3591 (electronic), 2001.

\bibitem{Tsuji}
M.~Tsuji.
\newblock {\em Potential theory in modern function theory}.
\newblock Maruzen Co. Ltd., Tokyo, 1959.

\bibitem{Waddington}
Simon Waddington.
\newblock Large deviation asymptotics for {A}nosov flows.
\newblock {\em Ann. Inst. H. Poincar\'e Anal. Non Lin\'eaire}, 13(4):445--484,
  1996.

\end{thebibliography}
\end{document}